\documentclass[reqno]{amsart}
\usepackage{amsmath}
\newtheorem{thm}{Theorem}
\newtheorem{prop}[thm]{Proposition}

\newtheorem{cor}[thm]{Corollary}
\newcommand{\R}{{\mathbb R}}

\newcommand{\C}{{\mathbb C}}
\title{A Note on Isometric Embeddings 
of Surfaces of Revolution \\}
\author{Martin Engman}
\thanks{Partially supported by the NSF Grant: Model Institutes for Excellence at 
UMET}
\date{}
\begin{document}
\address{Departamento de Ciencias y Tecnolog\'{i}a, Universidad Metropolitana, 
San Juan, PR 00928}
\email{um\_mengman@suagm.edu}

\maketitle

\noindent
{\bf 1. INTRODUCTION.} In 1916 Weyl asked: ``Does there exist a smooth global isometric embedding
 of the Riemannian manifold 
$(S^2,g)$ into $(\R^3, can)$?" Here, $can$ is the flat (i.e. canonical) metric on $\R^3$.
Solutions to this problem (of various orders of differentiability) have been given by,
 Weyl, \cite{we}, 
Nirenberg, \cite{ni}, Heinz, \cite{he}, Alexandrov, 
\cite{al}, and Pogorelov, \cite{po}. All of these results require that the Gauss curvature be 
positive.
A recent result of Hong and Zuily \cite{hz} requires only that the curvature be 
non-negative. The converses of these theorems cannot hold since it is easy to find embedded 
metrics with some negative curvature. Moreover, not every metric on $S^2$ admits such
an isometric embedding. The reader may refer to Greene, \cite{gr}, 
wherein one finds examples of smooth metrics on 
$S^2$ for which there is no $C^2$ isometric embedding in $(\R^3, can)$.
In this paper we study the special case of $S^1$ invariant metrics on $S^2$ 
(i.e. surfaces of revolution). These metrics have two important characteristics: they have two 
distinguished points, the 
fixed points of the $S^1$ action (north and south pole); and important geometric quantities, 
such as Gauss curvature, 
are functions of only one variable. This has the effect of converting computations, which in 
general require Stokes theorem, to the more basic fundamental theorem of calculus. 

Armed with nothing more than the fundamental theorem of 
calculus, we will show that the non-negativity of the
integrals of the curvature over all disks centered at a pole is the necessary and 
sufficient condition for isometric 
embeddability of $S^1$ invariant metrics on $S^2$ into $(\R^3, can)$. So that, although in general, an isometrically
embedded metric cannot be expected to have non-negative curvature everywhere, in some sense one 
has the next best thing: every pole-centered disk has non-negative integral curvature.
  
\vspace{.2in}
 
\noindent
{\bf2. ISOMETRIC EMBEDDINGS OF ROTATION INVARIANT SURFACES.} 
The results that we obtain are based on a well know and quite elementary result regarding 
$S^1$ invariant metrics on $S^2$ (i.e. Surfaces of Revolution)
which can be found in Besse \cite{be}, p. 95-106. The reader will find the results of this 
section to be, 
essentially, nothing but a reformulation of Besse's treament. The reader should be aware that 
our use of the term {\em surface of revolution} is equivalent to the definition {\em abstract 
surface of revolution} (not, {\em a priori}, embedded in $\R^3$) which can be found in Hwang, 
\cite{hwa}.

Let $(M,g)$ be an $S^1$ invariant Riemannian manifold which is diffeomorphic
to $S^2$ and which has volume $4\pi$. We will assume the metric to be
$C^{\infty}$.
 Since $(M,g)$ has an effective $S^{1}$ isometry group there are exactly two 
fixed points. We call the fixed points $np$ and $sp$ and let $U$ be the
open set $M \setminus \{np,sp\}$. On $U$ the metric has the form
\begin{equation}
 g= ds \otimes ds + a^{2}(s)d\theta \otimes d\theta   
\label{eq:feq}
\end{equation} 
where $s$ is the 
arclength along a geodesic connecting $np$ to $sp$ and $a(s)$ is a 
function $a:[0,L]\rightarrow \R^+$ satisfying $a(0)=a(L)=0$ and 
$a'(0)=1=-a'(L)$. The canonical (i.e. constant curvature) metric is obtained by taking 
$a(s)= \sin s$. 

The Riemannian measure, $ \sqrt{g} ds _{ \wedge} d \theta \stackrel{\rm def}{=} \sqrt{\det g} ds
_{ \wedge} d \theta$ (with
respect to which all surface integrals are computed as double integrals) is, in this case, given
by $\sqrt{g} ds_{ \wedge} d \theta=a(s) ds_{ \wedge} d \theta$. This fact, together with the 
following formula
for the Gauss curvature of this metric,
\begin{equation}
K(s) = -a''(s)/a(s),
\end{equation}
reduces the proof of the Gauss-Bonnet Theorem, $\int_M K = 4\pi$, to a simple application of the 
fundamental theorem of calculus.
 
One can always isometrically embed such a metric into $\R^2 \times \C$  as follows:
\begin{equation}
\left\{ \begin{array}{ccl}
  \psi^1(s,\theta) & = & a(s)\cos \theta \\
  \psi^2(s,\theta) & = & a(s)\sin \theta \\
  \psi^3(s,\theta) & = & \int_c^s \sqrt{1-(a')^2(t)}dt
                      \end{array} \right. \label{eq:param}
\end{equation} 
for some $c \in [0,L]$. The function $\psi^3$ is, in general, complex valued so the 
condition for embeddability in $\R^3$ reduces to the condition under which $\psi^3$ is 
real valued. 
The interested reader will find the following proposition to be 
equivalent to the result on p-106 of \cite{be}, so we will omit the proof.

\begin{prop} Let $(M,g)$, with metric $g$ as in \eqref{eq:feq}, be 
diffeomorphic to $S^2$. $(M,g)$ can be isometrically $C^1$ embedded 
in $(\R^3,can)$ if and only if $|a'(s)| \leq 1$ for all $s \in [0,L]$. 
\qed
\end{prop}
For a generalization of this proposition, see \cite{hwa}.

We will end this section with the observation
that, in our special case, non-negative curvature implies isometric 
embeddability since $K(s) \geq 0$ implies that $a'(s)$ is a non-increasing
function on $[0,L]$ with maximum $1$ and minimum $-1$ and hence $|a'(s)| \leq 1$. This proves 
the following proposition which can be though of as a special case of the theorem of 
Hong and Zuily \cite{hz}.

\begin{prop} Let $(M,g)$ be an $S^1$ invariant Riemannian manifold which is
diffeomorphic to $S^2$ with Gauss curvature $K$. If $K \geq 0$ on $M$ then $(M,g)$ can be $C^1$ 
isometrically embedded in $(\R^3,can)$.
\qed
\end{prop}
\vspace{.2in}

\noindent
{\bf 3. MAIN RESULTS.} Our main result is the following necessary and sufficient condition. 

\begin{thm} 
Let $(M,g)$ be an $S^1$ invariant Riemannian manifold which is
diffeomorphic to $S^2$. Let $K$ be the Gauss curvature and let $np$ and $sp$ denote the fixed 
points of the $S^1$ action. $(M,g)$ can be $C^1$ isometrically embedded in $(\R^3,can)$ if and only 
if for all disks $\Omega$  centered at $np$ or $sp$  $$\int_{\Omega} K \geq 0.$$ 

\end{thm}

\begin{proof} 
Each closed disk, $\Omega$, centered at $np$ is parametrized by $0 \leq s \leq x$, 
$0 \leq \theta < 2 \pi$ for 
some $x$. So 
$$ \int_{\Omega} K = \int^{2\pi}_0 \int^x_{0} K(s) \sqrt{g} ds_{ \wedge} d \theta =2\pi 
\int^x_{0} K(s) \sqrt{g} ds.$$
Hence 
$$\int_{\Omega} K  \geq 0 \mbox{ if and only if } \int^x_{0} K(s) \sqrt{g} ds \geq 0.$$
But now, since $K(s) = -a''(s)/a(s)$, $\sqrt{g}=a(s)$, and 
$a'(0)=1$,
\begin{equation}
 \int^x_{0} K(s)\sqrt{g} ds= -a'(x)+1 \label{eq:K}
\end{equation}
and it now easily follows that 
$$\int^x_{0} K(s) \sqrt{g} ds \geq 0 \mbox{ if and only if } a'(x) \leq 1.$$

A similar argument shows that the non-negativity of all integrals over disks centered at $sp$ is 
equivalent to 
$\int^L_{x} K(s) \sqrt{g} ds \geq 0$ for all $x$. Also, it is easy to show that $\int^L_{0} K(s) \sqrt{g} ds 
=2$ (this is, essentially, the Gauss-Bonnet
Theorem) hence 
$$\int^L_{x} K(s) \sqrt{g} ds = 2-\int^x_{0} K(s) \sqrt{g} ds$$
which, after substitution of \eqref{eq:K}, leads to the conclusion 
$$\int^L_{x} K(s) \sqrt{g} ds \geq 0 \mbox{ if and only if } a'(x) \geq -1.$$
In other words: $\int_{\Omega} K \geq 0$ for all disks $\Omega$  centered at $np$ or $sp$ if and only if 
$|a'(s)| \leq 1$ for all $s \in [0,L]$ and the result now follows from Proposition 1.
\end{proof}

There are several interesting consequences of this result. Perhaps the first observation to make 
is that Theorem 3 provides another proof of Proposition 2 since, obviously, $
K \geq 0$ everywhere implies $\int_{\Omega} K \geq 0$ for all $\Omega$.  
It is also easy to see that if $(M,g)$ is isometrically embedded in $(\R^3,can)$, then the 
curvature must be non-negative at the poles.
This observation yields the following obstruction to isometric embedding in terms of the 
curvature. 
 For a 
similar result in terms of the spectrum see \cite{eng}.

\begin{cor} 

Let $(M,g)$ be an $S^1$ invariant Riemannian manifold of area $4\pi$ which is
diffeomorphic to $S^2$. Let $K$ be the Gauss curvature and let $np$ and $sp$ denote the fixed points of the $S^1$ action. 
If $K(np)<0$ or $K(sp)<0$ then $(M,g)$ cannot be isometrically embedded in $(\R^3,can)$. 

\end{cor}

\begin{proof}
We assume $K(np)<0$ i.e. $K(0)<0$. Since $K(s)$ extends to a continuous function on $[0,L]$ and 
$\sqrt{g(s)}>0$ on $(0,L)$ then for some 
$\epsilon > 0$, 
$K(s)\sqrt{g(s)}<0$ on $(0,\epsilon)$ and hence $\int^{\epsilon}_{0} K(s)\sqrt{g(s)} ds < 0$. The
result now follows from Theorem 3.
\end{proof}

This result is probably well known, but we have not found a reference.

Another particular case of the Gauss Bonnet theorem states that for any domain $\Omega$ with 
$C^2$ boundary $C$ 
\begin{equation}
\int_{\Omega} K + \oint_C k_g =2\pi \label{eq:L}
\end{equation}
where $k_g$ is the geodesic curvature of $C$. Applying this to our disks centered at the pole 
and observing that the boundary of these disks are circles of latitude, we
get:

\begin{cor}

Let $(M,g)$ be an $S^1$ invariant Riemannian manifold of area $4\pi$ which is
diffeomorphic to $S^2$ and has Gauss curvature $K$. $(M,g)$ can be $C^1$ isometrically embedded 
in $(\R^3,can)$ if and 
only if for all circles of latitude, $C$, $$\left| \oint_C k_g \right| \leq 2\pi.$$ 

\end{cor}

\begin{proof} Each circle of latitude serves as the boundary of two disks, $\Omega_1$ centered at
 $np$ and another, 
$\Omega_2$ centered at $sp$ but with opposite orientations hence, from \eqref{eq:L}, 

$$ \int_{\Omega_1} K + \oint_C k_g =2\pi \mbox{ and } \int_{\Omega_2} K - \oint_C k_g =2\pi,$$
and so 
$\int_{\Omega} K \geq 0$ for all $\Omega$ if and only if $2\pi - \oint_C k_g 
\geq 0$ and $2\pi + \oint_C k_g 
\geq 0$. Thus $\left| \oint_C k_g \right| \leq 2\pi$.
\end{proof}

\end{document}